\documentclass{amsart}
\usepackage{graphicx}
\usepackage{amssymb}
\usepackage{amsmath}
\usepackage{amsthm}
\usepackage{subfigure}
\usepackage{enumerate}
\usepackage{xcolor}
\usepackage{bbm}
\usepackage{calligra}
\DeclareMathAlphabet{\mathcalligra}{T1}{calligra}{m}{n}

\newtheorem{thm}{Theorem}[section]
\newtheorem{cor}[thm]{Corollary}

\newtheorem{prop}[thm]{Proposition}
\newtheorem{prob}[thm]{Problem}
\theoremstyle{definition}
\newtheorem{defn}[thm]{Definition}
\newtheorem{rem}[thm]{Remark}
\newtheorem*{defn*}{Definition}
\newtheorem*{rems*}{Remarks}
\newtheorem*{rem*}{Remark}

\numberwithin{equation}{section}

\newcommand{\Eq}{{\text{E}}}

\newcommand{\M}{{M}}

\newcommand{\F}{\text{F}}

\newcommand{\Cwms}{\text{CWMS}}
\newcommand{\Sms}{\text{SMS}}
\renewcommand{\d}{\mathrm{d}}

\begin{document}

\title[CWMS, SMS and isoperimetric equalities] {The Constant Width Measure Set, \linebreak the Spherical Measure Set \linebreak and isoperimetric equalities for planar ovals}
\author{ Micha\l{} Zwierzy\'nski}
\address{Warsaw University of Technology\\
Faculty of Mathematics and Information Science\\
ul. Koszykowa 75\\
00-662 Warszawa\\
Poland}

\email{Michal.Zwierzynski@pw.edu.pl}
\thanks{The work of M. Zwierzy\'nski was partially supported by NCN grant no. DEC-2013/11/B/ST1/03080. }



\begin{abstract}

In this paper we introduce the \textit{Constant Width Measure Set}, which measures the constant width property of an oval, i.e. the planar simple closed strictly convex curve. We study its geometrical properties.

We find the exact relation between the length and the area of the region bounded by an oval $\M$. Namely, the following equality is fulfilled:
\begin{align*}
L_{\M}^2 &=4\pi A_M+8\pi\left|\widetilde{A}_{\Eq_{\frac{1}{2}}(\M)}\right|+\pi\left|\widetilde{A}_{\Cwms(\M)}\right|,
\end{align*}
where $L_{\M}, A_{\M}, \widetilde{A}_{\Eq_{\frac{1}{2}}(\M)}, \widetilde{A}_{\Cwms(\M)}$ are the length of $\M$, the area bounded by $\M$, the oriented area of the Wigner caustic of $\M$ and the oriented area of the Constant Width Measure Set of $\M$, respectively. 

Furthermore we study the geometry of the \textit{Spherical Measure Set}, which is an offset of a curve with a special distance. We show that the oriented area of this set of an oval $\M$, $\widetilde{A}_{\Sms(\M)}$, satisfies the following equality:
\begin{align*}
4\left|\widetilde{A}_{\Sms(\M)}\right|=8\left|\widetilde{A}_{\Eq_{\frac{1}{2}}(\M)}\right|+\left|\widetilde{A}_{\Cwms(\M)}\right|.
\end{align*}
\end{abstract}

\maketitle

\section{Introduction}

The \textit{classical isoperimetric problem} in the Euclidean plane $\mathbb{R}^2$ states that:

\begin{prob}(Isoperimetric problem)
Find the figure which has the maximum area for a given perimeter.
\end{prob}

The solution (without any proof) was known by the ancient Greeks. The first mathematical proof that the disc is the unique answer to the isoperimetric problem was given by Steiner in  \cite{S1}. He proved the isoperimetric inequality.

\begin{thm}(Isoperimetric inequality)
Let $\M$ be a simple closed curve of the length $L_{\M}$, enclosing a region of the area $A_{\M}$, then
\begin{align}\label{IsoperimetricIneq}
L_{\M}^2\geqslant 4\pi A_{\M},
\end{align}

and the equality (\ref{IsoperimetricIneq}) holds if and only if $\M$ is a circle.

\end{thm}

In \cite{Z2} we proved the improved isoperimetric inequality, which also gives the isoperimetric equality for convex curves of constant width.

\begin{thm}(Improved isoperimetric inequality 1)\label{ImprovedIsoperimetricIneq1}
Let $\M$ be a closed regular simple convex curve. Then
\begin{align*}
L_{\M}^2\geqslant 4\pi A_{\M}+8\pi\left|\widetilde{A}_{E_{\frac{1}{2}}(\M)}\right|,
\end{align*}
where $\widetilde{A}_{E_{\frac{1}{2}}(\M)}$ is the oriented area of the Wigner caustic of $\M$, and the equality holds if and only if $\M$ is a curve of constant width.
\end{thm}

An affine equidistant is the set of points of chords connecting points on $\M$ where tangent lines to $\M$ are parallel, which divide the chord segments between the base points with a fixed ratio $\lambda$. If the ratio $\lambda$ is equal to $\displaystyle\frac{1}{2}$ then this set is also known as the \textit{Wigner caustic} of $\M$. The Wigner caustic was first introduced by Berry in \cite{B1}. There are many papers considering affine equidistants and in particular the Wigner caustic, see \cite{C2, DJRR1, DMR1, DR1, DRS1, DRZsecant, DZ1, DZnew, DZ2, DZgb, G3, GReeve, GWZ1, JJR1, RZ1, Z1, Z2}, and the literature therein. The Wigner caustic is also known as the \textit{area evolute} (see, e.g. \cite{C2, G3}), as the \textit{symmetry defect} (for example, see \cite{DFJSymmetryDefect, JJR1}) and as the \textit{middle hedgehog} (\cite{S5, S6}). Singularities of the Wigner caustic for ovals occur exactly from an antipodal pair (the tangents at the two points are parallel and the curvatures are equal), the well-known Blaschke-S\"uss theorem states that there are at least three pairs of antipodal points on an oval (\cite{G1, L2}). 

Offset curves and surfaces are well-known geometric objects in the field of mathematics and computer aided geometric design, possibly because they give a powerful tool in many applications (see \cite{FHK1, HL1, KGP1, PP1, Rochera2}, and the references therein). The Spherical Measure Set is an offset.

The classical four-vertex theorem states that the curvature function of a simple, closed, smooth plane curve has at least four local extrema. Four-Vertex Theorem was first proved for ovals in 1909 by Mukhopadhyaya (\cite{M1}) and in general by Kneser in 1912 using a projective argument (\cite{K1}), see also \cite{DGP1, NBRF1} and the literature therein. We present a new proof of Four-Vertex Theorem for generic ovals and the proof is a simply consequence of the number of singular points of the Spherical Measure Set of an oval -- see Corollary \ref{Cor4VertexThm}.

There are many important inequalities in the convex geometry and differential geometry, such as the isoperimetric inequality, the Brunn-Minkowski inequality, Aleksandrov-Fenchel inequality, Gage's inequality. The stability property of them are of great interest in geometric analysis. One of the methods used is the method of Fourier series and spherical harmonics based on the support function. It is worth mentioning that in recent years, many researchers have been focusing on works related to support functions. For the references, see see \cite{CGR1, F2, G6, G2, G4, H3, PX1, Rochera1, S2, YY, ZD} and the literature therein.

The paper is organized as follows.

In Section 2 we present geometric quantities, their Fourier series and affine equidistants, including the Wigner caustic.

Section 3 contains a definition of the Constant Width Measure Set and properties of this set for planar ovals.

In Section 4 we introduce the Spherical Measure Set and present properties of this set for planar simple regular closed curves.

Section 5 is devoted to prove the main results, isoperimetric equalities for planar simple regular closed curves.

In Section 6 we study the stability property of the isoperimetric inequality involving the Constant Width Measure Set and thanks to it we find the lower bounds of the absolute value of the oriented area of the Wigner caustic, the Constant Width Measure Set and the Spherical Measure Set.


\section{Geometric quantities, affine equidistants and the Fourier series}

Let $\M$ be a smooth planar curve, i.e. the image of the $C^{\infty}$ smooth map from an interval to $\mathbb{R}^2$. A smooth curve is called \textit{closed} if it is the image of a $C^{\infty}$ smooth map from $S^1$ to $\mathbb{R}^2$. A smooth curve is called \textit{regular} if its velocity does not vanish. A regular closed curve is called \textit{convex} if its signed curvature has a constant sign. An \textit{oval} is a smooth regular closed convex curve which is simple, i.e. it has no self-intersections.

\begin{defn}\label{parallelpair}
Let $\M$ be a regular closed curve. A pair $a,b\in\M$ ($a\neq b$) is called the \textit{parallel pair} if tangent lines to $\M$ at points $a,b$ are parallel.
\end{defn}

\begin{defn}\label{equidistantSet}
For a given $\lambda\in\mathbb{R}$, an affine $\lambda$-equidistant is the following set:
$$\Eq_{\lambda}(\M)=\left\{\lambda a+(1-\lambda)b\ \big|\ a,b \text{ is a parallel pair of } \M\right\}.$$

The set $\Eq_{\frac{1}{2}}(\M)$ is known as the \textit{Wigner caustic} of $\M$.
\end{defn}

\begin{figure}[h]
\centering
\includegraphics[scale=0.22]{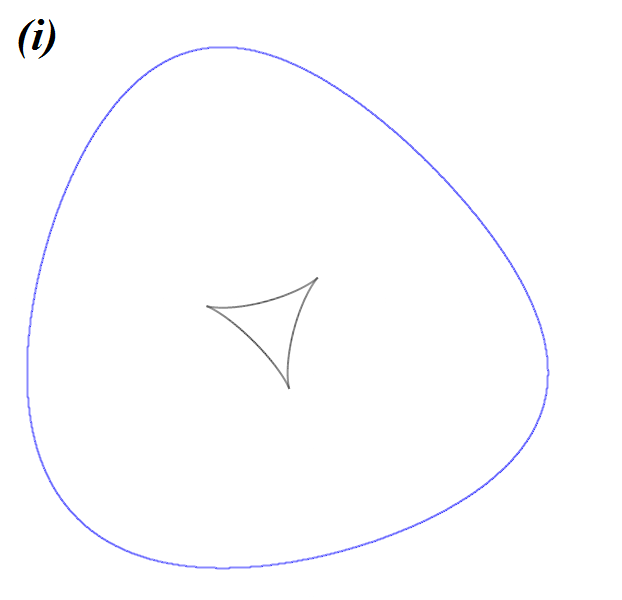}
\includegraphics[scale=0.22]{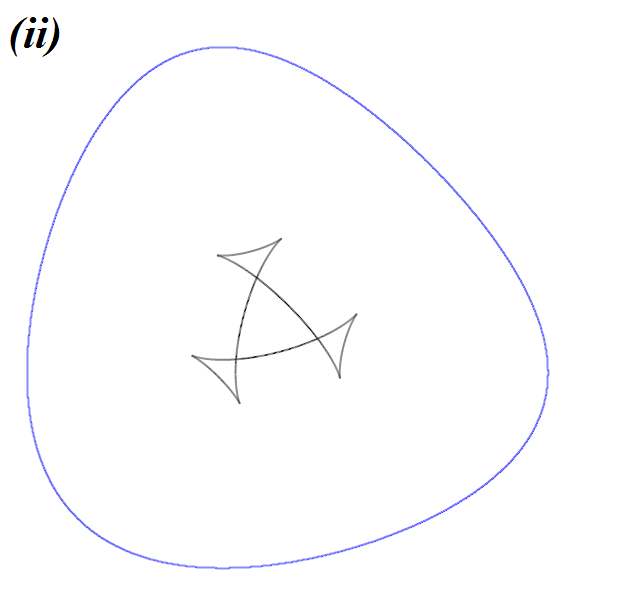}
\caption{An oval $\M$ and (i) $\Eq_{\frac{1}{2}}(\M)$, (ii) $\Eq_{\frac{2}{5}}(\M)$}
\label{Picture1}
\end{figure}

\begin{rem}\label{RemKnownWcEq}
It is known that if $\M$ is a generic oval, then for a generic $\lambda$ the set $\Eq_{\lambda}(\M)$ is a smooth closed curve with at most cusp singularities \cite{B1, G3}, the number of cusps of $\Eq_{\frac{1}{2}}(\M)$ is odd and not smaller than $3$ \cite{B1, G3} and the number of cusps of $\Eq_{\lambda}(\M)$ for a generic value of $\displaystyle\lambda\neq\frac{1}{2}$ is even \cite{DZ1}.
\end{rem}

\begin{defn}
An oval $\M$ is said to be a \textit{curve of constant width} if the distance between every pair of parallel tangent lines to $\M$ is constant. This constant is called the \textit{width} of the curve and we denote it by $w_{\M}$.
\end{defn}

\begin{defn}
An \textit{average width} of an oval $\M$ is $\displaystyle\overline{w}_{\M}=\frac{L_{\M}}{\pi}$.
\end{defn}

If $M$ is an oval of constant width, then by Barbier's theorem (\cite{B2}) $\overline{w}_{\M}=w_{\M}$.

The following details on plane ovals can be found in the classical literature \cite{G4}.

Let $\M$ be a positively oriented oval. Take a point $\Theta$ inside $\M$ as the origin of our frame. Let $p$ be the oriented perpendicular distance from $\Theta$ to the tangent line at a point on $\M$, and $\theta$ the oriented angle from the positive $x_1$-axis to this perpendicular ray. Clearly, $p$ is a single-valued periodic function of $\theta$ with period $2\pi$ and the parameterization of $\M$ in terms of $\theta$ and $p(\theta)$ is as follows
\begin{align}\label{ParameterizationM}
\gamma(\theta)=\big(\gamma_1(\theta),\gamma_2(\theta)\big)=\big(p(\theta)\cos\theta-p'(\theta)\sin\theta, p(\theta)\sin\theta+p'(\theta)\cos\theta\big).
\end{align}
The couple $\big(\theta, p(\theta)\big)$ is usually called the \textit{polar tangential coordinate} on $\M$, and $p(\theta)$ its \textit{Minkowski support function}.

Then, the curvature $\kappa$ of $\M$ is in the following form:
\begin{align}\label{CurvatureM}
\displaystyle \kappa(\theta)=\frac{1}{\rho(\theta)}=\frac{\d\theta}{\d s}=\frac{1}{p(\theta)+p''(\theta)}>0,
\end{align}
where $\rho$ denotes the radius of curvature of $\M$.

Let $L_{\M}$ and $A_{\M}$ be the length of $\M$ and the area it bounds, respectively. Then 
\begin{align}\label{CauchyFormula}
L_{\M}=\int_{\M}\d s=\int_0^{2\pi}\rho(\theta)\d\theta=\int_0^{2\pi}p(\theta)\d\theta,
\end{align}
and
\begin{align}\label{BlaschkeFormula}
A_{\M} & =\frac{1}{2}\int_{\M}\rho(\theta)\d s=\frac{1}{2}\int_0^{2\pi}p(\theta)\left[p(\theta)+p''(\theta)\right]\d\theta=\frac{1}{2}\int_0^{2\pi}\left[p^2(\theta)-p'^2(\theta)\right]\d\theta.
\end{align}
The formula (\ref{CauchyFormula}) (respectively (\ref{BlaschkeFormula})) is known as the \textit{Cauchy formula} (respectively the \textit{Blaschke formula}).

Since the Minkowski support function of $\M$ is continuous, bounded and $2\pi$-periodic, its Fourier series is in the form:
\begin{align}\label{Fourierofp}
p(\theta)=a_0+\sum_{n=1}^{\infty}\big(a_n\cos n\theta+b_n\sin n\theta\big).
\end{align}
Hence
\begin{align}\label{Fourierofpprime}
p'(\theta)=\sum_{n=1}^{\infty}n\big(-a_n\sin n\theta+b_n\cos n\theta\big).
\end{align}
By (\ref{Fourierofp}), (\ref{Fourierofpprime}) we get that
\begin{align}
\label{Lengthofmfourier} L_{\M} &=2\pi a_0,\\
\label{Areaofmfourier} A_{\M} &=\pi a_0^2-\frac{\pi}{2}\sum_{n=2}^{\infty}(n^2-1)(a_n^2+b_n^2).
\end{align}

Since $\gamma(\theta), \gamma(\theta+\pi)$ is a parallel pair of $\M$, $\gamma_{\lambda}$ -- the parameterization of $\Eq_{\lambda}(\M)$, is as follows:
\begin{align}\label{ParameterizationEqM}
\gamma_{\lambda}(\theta) 
	&=\big(\gamma_{\lambda, 1}(\theta), \gamma_{\lambda, 2}(\theta)\big)=\lambda\gamma(\theta)+(1-\lambda)\gamma(\theta+\pi) \\
\nonumber	&=\left(P_{\lambda}(\theta)\cos\theta-P'_{\lambda}\sin\theta, P_{\lambda}(\theta)\sin\theta+P'_{\lambda}(\theta)\cos\theta\right),
\end{align}
where $P_{\lambda}(\theta)=\lambda p(\theta)-(1-\lambda)p(\theta+\pi)$, $\theta\in[0,2\pi]$. 

\begin{rem}\label{RemWcDouble} If $\displaystyle\lambda=\frac{1}{2}$, then the map $\M\ni\gamma(\theta)\mapsto\gamma_{\frac{1}{2}}(\theta)\in\Eq_{\frac{1}{2}}(\M)$ for $\theta\in[0,2\pi]$ is the double covering of the Wigner caustic of $\M$.
\end{rem}

In particular
\begin{align}\label{SupportWc}
P_{\frac{1}{2}}(\theta) &= \frac{1}{2}(p(\theta)-p(\theta+\pi))\\
\label{SupportFourierWc} &= \sum_{\substack{n=1,\\ n\text{ is odd}}}^{\infty}\big(a_n\cos n\theta+b_n\sin n\theta).
\end{align}

\begin{defn}
Let $C$ be an oriented closed curve. Then the \textit{oriented area} of $C$ is
\begin{align*}
\widetilde{A}_C:=\frac{1}{2}\int_C-ydx+xdy=\iint_{\mathbb{R}^2}w_C(x,y)dxdy,
\end{align*}
where $w_C(x,y)$ is the winding number of $C$ around a point $(x,y)\in\mathbb{R}^2$.
\end{defn}

One can check that the oriented area of the Wigner caustic has the following formula in terms of coefficients of the Fourier series of the Minkowski support function $p$.

\begin{align}
\label{WCAreaFourier}2\widetilde{A}_{\Eq_{\frac{1}{2}}(\M)} &=\frac{1}{2}\int_0^{2\pi}\Big[P^2_{\frac{1}{2}}(\theta)-P_{\frac{1}{2}}'^2(\theta)\Big]\d\theta\\
\nonumber 
	&=-\frac{\pi}{2}\sum_{\substack{n=2,\\ n\text{ is odd}}}^{\infty}(n^2-1)(a_n^2+b_n^2).
\end{align}

\begin{rem}
Let us notice that $\M$ is an oval of constant width if and only if coefficients $a_{2n}, b_{2n}$ for $n\geqslant 1$ in the Fourier series of the Minkowski support function $p$ are all equal to zero \cite{F1, G4}.
\end{rem}

Wigner caustic of an oval $\M$ is an example of a planar \textit{hedgehog}, i.e. a planar curve which can be viewed as a Minkowski's difference of convex planar bodies (see \cite{MMY}).


\section{The Constant Width Measure Set}
\begin{defn}\label{DefCWMS}
Let $\M$ be a positively oriented oval. The \textit{Constant Width Measure Set} of $\M$ is the following set:
\begin{align}\label{CWMSDef}
\Cwms(\M)=\Big\{a-b+\overline{w}_{\M}\cdot\mathbbm{n}(a)\ \Big|\ \ a,b\text{ is a parallel pair of }\M\Big\},
\end{align}
where $\overline{w}_{\M}$ is a average width of $\M$ and $\mathbbm{n}(a)$ is a continuous unit normal vector field to $\M$ at $a$ compatible with the orientation of $\M$. We treat $\Cwms(\M)$ as a subset of a vector space $V=\mathbb{R}^2$. Let us assume that $\Theta$ is the origin of $V$. 
\end{defn}

\begin{figure}[h]
\centering
\includegraphics[scale=0.20]{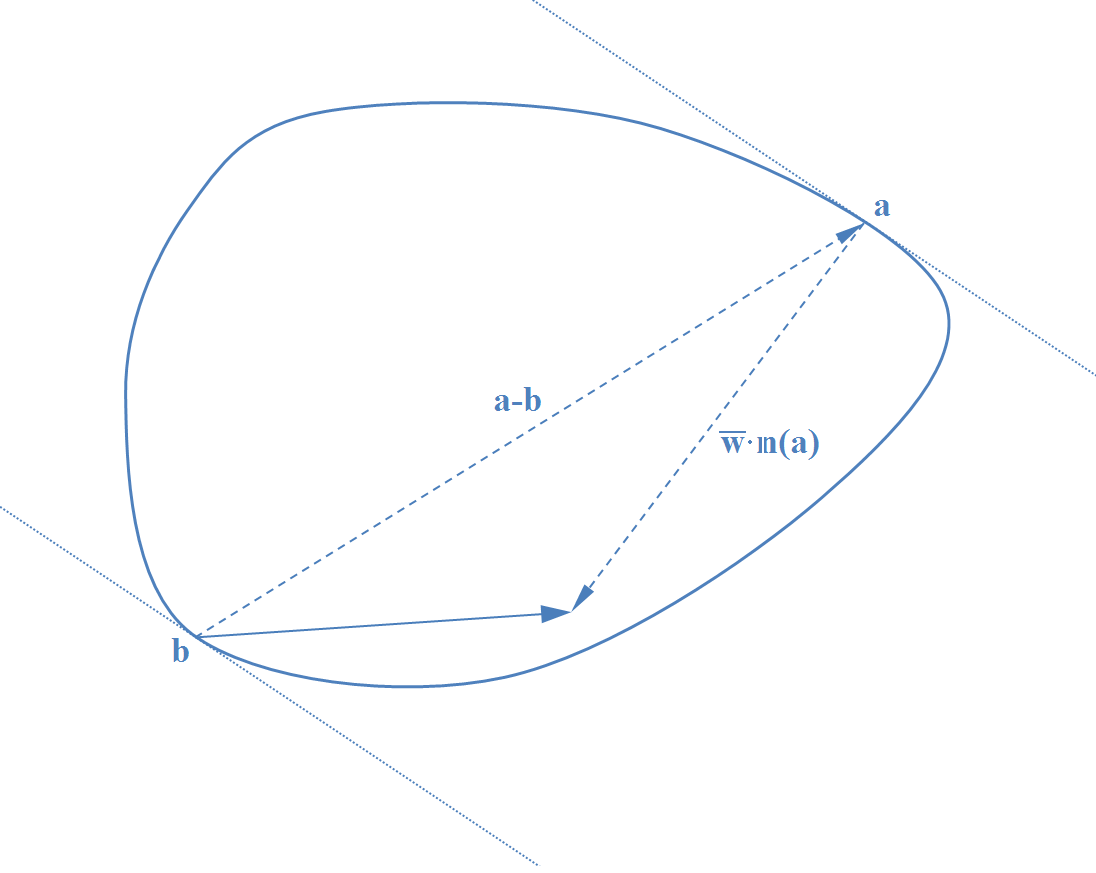}
\caption{An oval $\M$, a parallel pair $a,b$, tangent lines to $\M$ at $a$ and $b$, vectors $a-b, \overline{w}_{\M}\cdot\mathbbm{n}(a)$ and $a-b+\overline{w}_{\M}\cdot\mathbbm{n}(a)$}
\label{PictureCWMSDefn}
\end{figure}

Based on the preprint of our work and our Definition \ref{DefCWMS}, other researchers have introduced sets analogous to the $\Cwms$ set -- specifically for rosettes (see \cite{Z3}) and in normed spaces (see \cite{SC1}). They have also studied their properties and provided applications in their respective works.

\begin{figure}[h]
\centering
\includegraphics[scale=0.25]{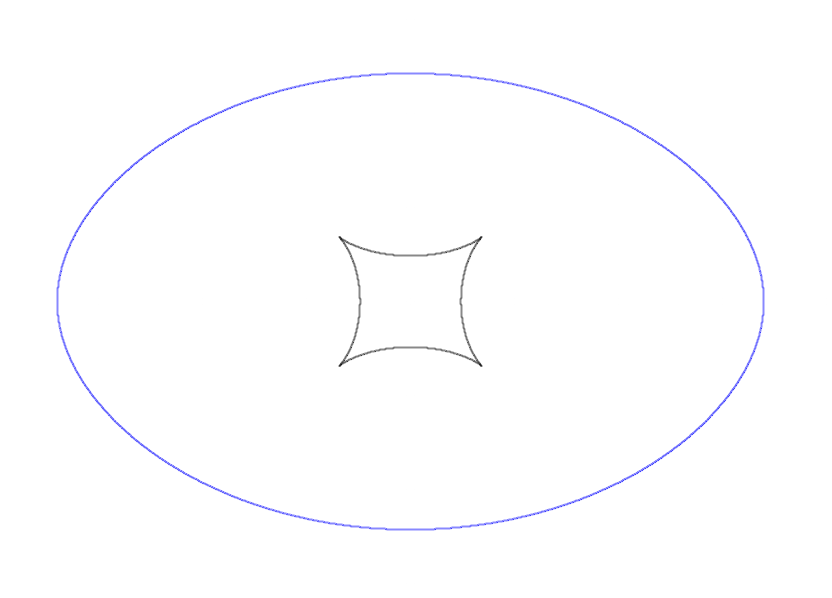}
\caption{An oval $\M$ and $\Cwms(\M)$}
\label{PictureCWMSEx}
\end{figure}

In Figure \ref{PictureCWMSDefn} we illustrate an example of a vector which belongs to $\Cwms(\M)$ and Figure \ref{PictureCWMSEx} illustrates $\Cwms(\M)$ of an oval $\M$.

\begin{thm}
$\M$ is an oval of constant width if and only if $\Cwms(\M)=\{\Theta\}$.
\end{thm}
\begin{proof}
If $\Cwms(\M)=\{\Theta\}$, then $b-a=\overline{w}_{\M}\cdot\mathbbm{n}(a)$ for all parallel pairs $a,b$. This gives that the length of $b-a$ is always constant and equal to $\overline{w}_{\M}$ and consequently $\M$ is an oval of constant width.

If $\M$ is an oval of constant width $w_{\M}$, then one can notice that if $a,b$ is a parallel pair of $\M$, then $a-b$ is perpendicular to the tangent lines to $\M$ at $a$ and $b$. This gives $a-b+w_{\M}\cdot\mathbbm{n}(a)=\Theta$.
\end{proof}

\begin{thm}\label{ThmCurvParallelCwms}
Let $\M$ be a positively oriented oval. Let $\kappa_{\M}(a)$ denote the curvature of $\M$ at $a$. Let $a,b$ be a parallel pair of $\M$ and let $q=a-b+\overline{w}_{\M}\cdot\mathbbm{n}(a)$ be a non-singular point of $\Cwms(\M)$. Then 
\begin{enumerate}[(i)]
\item the tangent line to $\Cwms(\M)$ at $q$ is parallel to the tangent lines to $\M$ at $a$ and $b$.
\item the curvature of $\Cwms(\M)$ at $q$ is equal
\begin{align}\label{KappaCmws}
\kappa_{\Cwms(\M)}(q)=\frac{\kappa_{\M}(a)\cdot\kappa_{\M}(b)}{|\kappa_{\M}(a)+\kappa_{\M}(b)-\overline{w}_{\M}\cdot\kappa_{\M}(a)\cdot\kappa_{\M}(b)|}.
\end{align}
\end{enumerate}
\end{thm}
\begin{proof}
Let (\ref{ParameterizationM}) be the parameterization of $\M$ and $p(\theta)$ its Minkowski support function.

One can check that the Minkowski support function of $\Cwms(\M)$ is in the following form:
\begin{align}\label{SupportCwms}
p_{\Cwms(\M)}(\theta)=p(\theta)+p(\theta+\pi)-\overline{w}_{\M}.
\end{align}
The formula (\ref{SupportCwms}) in terms of coefficients of a Fourier series of $p(\theta)$ is as follows:
\begin{align}\label{SupportFourierCwms}
p_{\Cwms(\M)}(\theta)=2\sum_{\substack{n=2,\\ n\text{ is even}}}^{\infty}(a_n\cos n\theta+b_n\sin n\theta).
\end{align}
Hence the parameterization of $\Cwms(\M)$ is in the following form:
\begin{align}
\gamma_{\Cwms(\M)}(\theta) = \Big( &p_{\Cwms(\M)}(\theta)\cos\theta-p'_{\Cwms(\M)}(\theta)\sin\theta, \\
\nonumber	& p_{\Cwms(\M)}(\theta)\sin\theta+p'_{\Cwms(\M)}(\theta)\cos\theta\Big).
\end{align}

Then: 

\begin{align*}
\gamma'(\theta) &=\Big(p(\theta)+p''(\theta)\Big)\cdot\Big(-\sin\theta,\cos\theta\Big),\\
\gamma'_{\Cwms(\M)}(\theta) &=\Big(p_{\Cwms(\M)}(\theta)+p''_{\Cwms(\M)}(\theta)\Big)\cdot\Big(-\sin\theta,\cos\theta\Big).
\end{align*}

Therefore the tangent line to $\M$ at $\gamma(\theta)$ is parallel to the tangent line to $\Cwms(\M)$ at non-singular point $\gamma_{\Cwms(\M)}(\theta)$.

By (\ref{Fourierofp}) one can get that
\begin{align}
\label{FourierRho}\rho_{\M}(\theta) &=p(\theta)+p''(\theta)=a_0-\sum_{n=1}^{\infty}(n^2-1)(a_n\cos n\theta+b_n\sin n\theta),\\
\label{FourierRhoPi}\rho_{\M}(\theta+\pi) 	&=a_0-\sum_{n=1}^{\infty}(-1)^n(n^2-1)(a_n\cos n\theta+b_n\sin n\theta).
\end{align}

By (\ref{SupportFourierCwms}), (\ref{FourierRho}), (\ref{FourierRhoPi}):
\begin{align}\label{RadiusOfCurvCwms}
\rho_{\Cwms(\M)}(\theta) &=|p_{\Cwms(\M)}(\theta)+p''_{\Cwms(\M)}|\\
\nonumber	&=\left|-2\sum_{\substack{n=2,\\ n\text{ is even}}}^{\infty}(n^2-1)(a_n\cos n\theta+b_n\sin n\theta)\right|\\
\nonumber &=\left|\rho_{\M}(\theta)+\rho_{\M}(\theta+\pi)-\overline{w}_{\M}\right|,
\end{align}
which ends the proof.
\end{proof}

By Theorem \ref{ThmCurvParallelCwms} we obtain the following corollaries.

\begin{cor}\label{CorConstantSign}
Let $\M$ be a positively oriented oval. Then the curvature of \linebreak $\Cwms(\M)$ is positive on each regular connected component of $\Cwms(\M)$.
\end{cor}

\begin{cor}\label{SingCor}
Let $\M$ be a positively oriented oval of an average width equal to $\overline{w}_{\M}$. Let $\mathbbm{n}$ be a continuous unit normal vector field to $\M$. Let $a,b$ in $\M$ be a parallel pair. Let $\rho_{\M}(a), \rho_{\M}(b)$ denote the radius of curvature of $\M$ at $a, b$, respectively. Then $\Cwms(\M)$ is singular at the point $a-b+\overline{w}_{\M}\cdot\mathbbm{n}(a)$ if and only if 
\begin{align}\label{SingCondition}
\rho_{\M}(a)+\rho_{\M}(b)=\overline{w}_{\M}.
\end{align}
\end{cor}

\begin{thm}\label{Thm4ncusps}
Let $\M$ be a generic oval. Then $\Cwms(\M)$ is a connected smooth curve with cusp singularities and the number of cusps is a positive multiple of $4$.
\end{thm}
\begin{proof}
$\Cwms(\M)$ is a connected smooth curve because of (\ref{SupportFourierCwms}). By Corollary \ref{SingCor} and by theory of Thom (1975) \cite{T1} one can get that generically $\Cwms(\M)$ has cusp singularity when (\ref{SingCondition}) holds.
By (\ref{RadiusOfCurvCwms}) the condition (\ref{SingCondition}) can be written in terms of the Fourier series of the Minkowski support function of $\M$ in the following way:
\begin{align}\label{CondFourierCusp}
\rho_{\Cwms(\M)}(\theta)=0\ \iff\ \sum_{\substack{n=2,\\ n\text{ is even}}}^{\infty}(n^2-1)(a_n\cos n\theta+b_n\sin n\theta)=0.
\end{align}
Therefore the number of cusps of $\Cwms(\M)$ is the number of zeros of (\ref{CondFourierCusp}) in the range $0$ to $2\pi$. Since $\rho_{\Cwms(\M)}(0)=\rho_{\Cwms(\M)}(\pi)$, the number of cusps in the range from $0$ to $\pi$ must be even. The function $\rho_{\Cwms(\M)}(\theta)$ is $\pi$-periodic, hence the number of cusps in the range from $0$ to $2\pi$ is a multiple of $4$. Moreover, the number of cusps must be at least $4$, since the Fourier series in (\ref{CondFourierCusp}) begins with terms in $\cos 2\theta$ and $\sin 2\theta$.
\end{proof}

\begin{defn}\label{DefNormal}
The \textit{tangent line to the} $\Cwms(\M)$ \textit{at a cusp point} $p$ is the limit of a sequence of $T_{q_n}\Cwms(\M)$ in $\mathbb{R}P^1$ for any sequence $q_n$ of regular points of $\Cwms(\M)$ converging to $p$.
\end{defn}
It is easy to see that this definition does not depend on the choice of a converging sequence of regular points. By Theorem \ref{ThmCurvParallelCwms}(i) we can see that the tangent line to $\Cwms(\M)$ at the cusp point $a-b+\overline{w}_{\M}\cdot\mathbbm{n}(a)$ is parallel to the tangent lines to $\M$ at $a$ and $b$.

It is easy to see that if $\M$ is a generic oval, then for any line $l$ in $\mathbb{R}^2$ there exists exactly two points $a,b\in\M$ in which tangent lines to $\M$ at $a$ and $b$ are parallel to $l$. Therefore we get the following corollary.

\begin{cor}\label{CorNumberOfTangentLines}
Let $\M$ be a generic oval. Then for any line $l$ in $\mathbb{R}^2$  there exists exactly two points $p, q\in\Cwms(\M)$ such that tangent lines to $\Cwms(\M)$ at $p$ and $q$ are parallel to $l$.
\end{cor}

We can define the normal vector to $\Cwms(\M)$ at any its point $a-b+\overline{w}_{\M}\cdot\mathbbm{n}(a)$ in the following way.

\begin{defn}
The normal vector to $\Cwms(\M)$ at $a-b+\overline{w}_{\M}\cdot\mathbbm{n}(a)$ is $\mathbbm{n}(a)$.
\end{defn}

Let us notice that the continuous normal vector field to $\Cwms(\M)$ at regular and cusp points is perpendicular to the tangent line to $\Cwms(\M)$. Using this fact and the above definition we define the rotation number in the following way.

\begin{defn}\label{RotationNumberCusps}
The \textit{rotation number} of a smooth closed curve $\M$ with at most cusp singularities is the number of rotation of its continuous normal vector field. 
\end{defn}

\begin{prop}\label{PropRotationCwms}
Let $\M$ be a generic oval. Then the absolute value of the rotation number of $\Cwms(\M)$ is equal to one.
\end{prop}
\begin{proof}
It is a consequence of Theorem \ref{Thm4ncusps} and Corollary \ref{CorNumberOfTangentLines}.
\end{proof}

\begin{prop}
If $\M$ is an oval, then $\Cwms(\M)$ has the center of symmetry.
\end{prop}
\begin{proof}
Let $a,b$ be a parallel pair of $\M$, then $p=a-b+\overline{w}_{\M}\cdot\mathbbm{n}(a)\in\Cwms(\M)$ and $q=b-a+\overline{w}_{\M}\cdot\mathbbm{n}(b)\in\Cwms(\M)$, and because $\mathbbm{n}(a)=-\mathbbm{n}(b)$ one can notice that $\Theta$ is the center of symmetry of $\Cwms(\M)$.
\end{proof}

\begin{prop}
There exists oval $\M$ for which the Constant Width Measure Set has exactly $4n$ cusp singularities.
\end{prop}
\begin{proof}
One can check that $p_{4n}(\theta)=4n^2+2+\cos 2n\theta$ is a Minkowski support function of an oval $\M_{4n}$ and $\Cwms(\M_{4n})$ has exactly $4n$ cusps.
\end{proof}

\begin{figure}[h]
\centering
\includegraphics[scale=0.35]{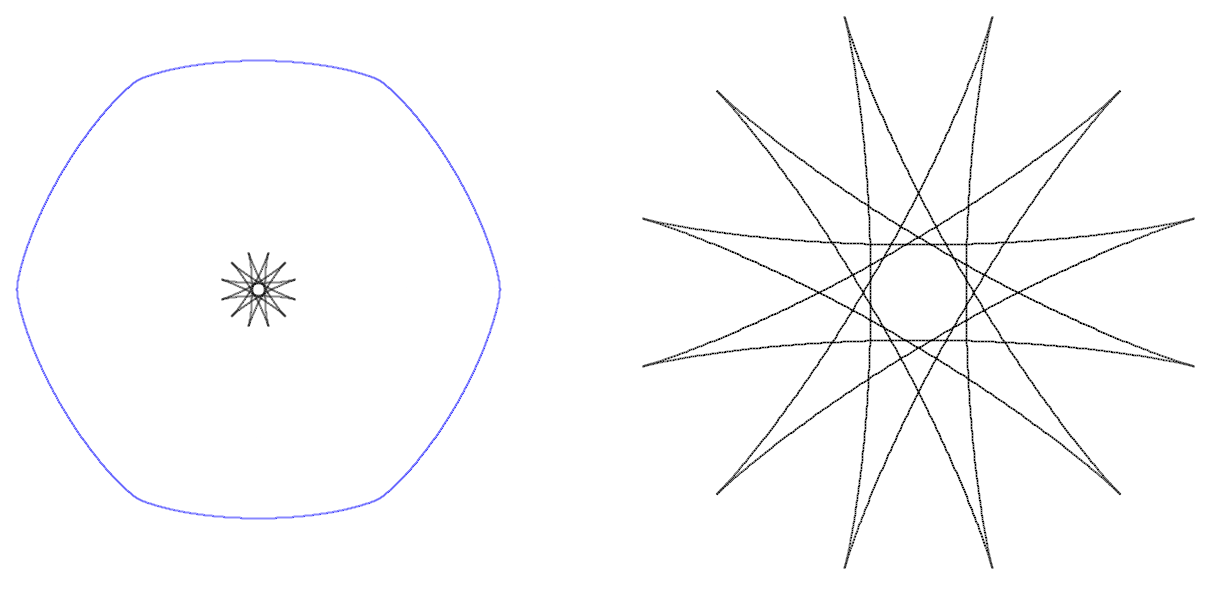}
\caption{An oval $\M_{12}$ and $\Cwms(\M_{12})$ which has $12$ cusp singularities. The support function of $M_{12}$ is $p_{12}(\theta)=38+\cos 6\theta$}
\label{PictureM12}
\end{figure}

In Figure \ref{PictureM12} we illustrate an example of an oval $\M$ for which the Constant Width Measure Set has exactly $12$ cusp singularities.

\begin{prop}
Let $\M$ be an oval. Then
\begin{align*}
L_{\Cwms(\M)}\leqslant 4L_{\M}.
\end{align*}
\end{prop}
\begin{proof}
Let $p(\theta)$ be the Minkowski support function of $\M$ and let $p_{\Cwms(\M)}$ (as in (\ref{SupportCwms})) be the support function of $\Cwms(\M)$.
Then
\begin{align*}
L_{\Cwms} &= \int_0^{2\pi}|\gamma'_{\Cwms(\M)}(\theta)|\d\theta\\
	&=\int_0^{2\pi}|p_{\Cwms(\M)}(\theta)+p''_{\Cwms(\M)}(\theta)|\d\theta\\
	&=\int_0^{2\pi}|p(\theta)+p''(\theta)+p(\theta+\pi)+p''(\theta+\pi)-\overline{w}|\d\theta\\
	&\leqslant\int_0^{2\pi}|p(\theta)+p''(\theta)|\d\theta+\int_0^{2\pi}|p(\theta+\pi)+p''(\theta+\pi)|\d\theta+\int_0^{2\pi}|\overline{w}|\d\theta\\
	&=L_{\M}+L_{\M}+2L_{\M}=4L_{\M}.
\end{align*}
\end{proof}

\section{The Spherical Measure Set}

\begin{defn}
Let $\M$ be a regular positively oriented simply closed curve. An \textit{$\alpha$-offset} for a given $\alpha\in\mathbb{R}$ of $\M$ is the following set:
\begin{align*}
\F_{\alpha}(\M)=\big\{a+\alpha\cdot\mathbbm{n}(a)\ \big|\ a\in\M\big\},
\end{align*}
where $\mathbbm{n}(a)$ is a continuous unit normal vector field to $\M$ at $a$ compatible with the orientation of $\M$.
\end{defn}

It is well known that offsets generically admit at most cusp singularities, singularities of all offsets of $\M$ form an evolute of $\M$ and set of all points of self-intersections of offsets forms a medial axis of $\M$. See \cite{FHK1, HL1, KGP1, PP1} and the literature therein.

\begin{defn}
The \textit{Spherical Measure Set} of a regular positively oriented simply closed curve $\M$ is an offset at level $\displaystyle \frac{L_{\M}}{2\pi}$,
\begin{align}
\Sms(\M)=\F_{\frac{L_{\M}}{2\pi}}(\M)=\left\{a+\frac{L_{\M}}{2\pi}\cdot\mathbbm{n}(a)\ \Big|\ a\in\M\right\},
\end{align}
where $\mathbbm{n}$ is the unit normal vector field compatible with the orientation of $\M$.
\end{defn}
\begin{rem}
If $\M$ is an oval, then $\displaystyle {\frac{L_{\M}}{2\pi}}={\frac{\overline{w}_{\M}}{2}}$.
\end{rem}

\begin{figure}[h]
\centering
\includegraphics[scale=0.55]{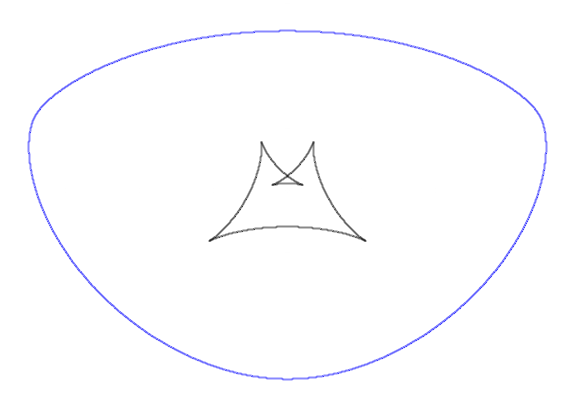}
\caption{An oval $\M$ and $\Sms(\M)$}
\label{PictureSms}
\end{figure}

\begin{figure}[h]
\centering
\includegraphics[scale=0.40]{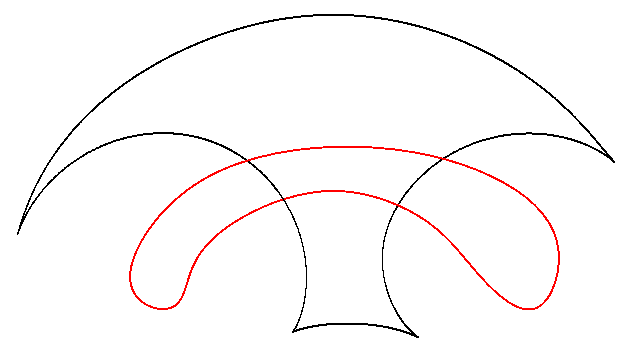}
\caption{A regular simple closed curve $\M$ and $\Sms(\M)$}
\label{PictureSms2}
\end{figure}

See Figure \ref{PictureSms} and Figure \ref{PictureSms2} for examples of $\Sms$.

\begin{prop}
Let $\M$ be a regular positively oriented simply closed curve. Then $\Sms(\M)=\{x\}$ if and only if $\M$ is a circle and $x$ is its center.
\end{prop}
\begin{proof}
If $\M$ is a circle and $x$ is its center, then it easy to see that $\frac{1}{2}\overline{w}$ is a radius of $\M$ and $\Sms(\M)=\{x\}.$

Let $\displaystyle w=\frac{L_{\M}}{\pi}$.
Let us assume that $\Sms(\M)=\{x\}$. Let $a,b$ be a parallel pair of $\M$, then $w\cdot\mathbbm{n}(a)=2x-2a$ and $w\cdot\mathbbm{n}(b)=2x-2b$. 
If $\mathbbm{n}(a)=\mathbbm{n}(b)$, then $a=b$, otherwise if $\mathbbm{n}(a)=-\mathbbm{n}(b)$, then $|a-b|=w$ and $a+b=2x$, so $\M$ is a curve of constant width and $x$ is the center of symmetry of $\M$, hence $\M$ is a circle.
\end{proof}

Bi direct calculations we get the following proposition.

\begin{prop}\label{ThmCurvParallelSms}
Let $\M$ be a positively oriented regular closed curve. Let $\kappa_{\M}$ denote the curvature of $\M$ at $a$. Then
\begin{enumerate}[(i)]
\item a point $\displaystyle a+\frac{L_{\M}}{2\pi}\cdot\mathbbm{n}(a)$ is a singular point of $\Sms$ if and only if 
\begin{align}\label{SingConditionSms}
L_{\M}\kappa_{\M}(a)=2\pi,
\end{align}
\item the tangent line to $\Sms(\M)$ at non-singular point $\displaystyle a+\frac{L_{\M}}{2\pi}\cdot\mathbbm{n}(a)$ is parallel to the tangent line to $\M$ at a point $a$.
\item the curvature of $\Sms(\M)$ at non-singular point $\displaystyle q=a+\frac{L_{\M}}{2\pi}\cdot\mathbbm{n}(a)$ is equal
\begin{align}\label{KappaSms}
\kappa_{\Sms(\M)}(q)=\frac{2\pi\kappa_{\M}(a)}{|2\pi - L_{\M}\kappa_{\M}(a)|}.
\end{align}
\end{enumerate}
\end{prop}

\begin{thm}\label{ThmCuspsSms}
Let $\M$ be a generic regular closed curve. Then the number of cusp singularities of $\Sms(\M)$ is even and not smaller than $2$.
If $\M$ is a generic oval, then the number of cusp singularities of $\Sms(\M)$ is not smaller than $4$.
\end{thm}

\begin{proof}

Let us assume that $\M$ is positively oriented.

\begin{figure}[h]
\centering
\includegraphics[scale=0.14]{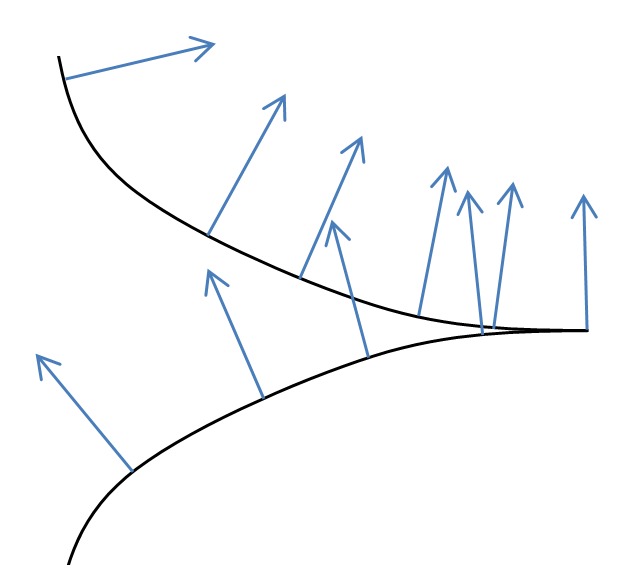}
\caption{The cusp singularity with a continuous normal vector field. Vectors in upper regular component of a curve are directed outside the cusp, others are directed inside the cusp}
\label{PictureNormalVectorToCusp}
\end{figure}

By  theory of Thom (1975) \cite{T1} one can get that generically $\Sms(\M)$ has cusp singularity when (\ref{SingConditionSms}) holds.

A continuous normal vector field to the germ of a curve with the cusp singularity is directed outside the cusp on the one of two connected regular components and is directed inside the cusp on the other component as it is showed in Figure \ref{PictureNormalVectorToCusp}.

Without loss of generality, let us assume that $\M$ is positively oriented. It is easy to see that the rotation number of $\Sms(\M)$ is equal to the rotation number of $\M$, which is an integer. Thus, the number of cusps of $\Sms(\M)$ is even.

Let us assume that $\Sms(\M)$ is regular. Let $[0,L_{\M})\ni s\mapsto\gamma(s)\in\mathbb{R}^2$ be the arc length parameterization of $\M$. Then by Proposition \ref{ThmCurvParallelSms}(i) we get that the curvature at each point of $\M$ is smaller or greather than $\displaystyle\frac{2\pi}{L_{\M}}$. If $\displaystyle\kappa_{\M}(s)<\frac{2\pi}{L_{\M}}$ for all $s$, then
\begin{align*}
\int_0^{L_{\M}}\kappa_{\M}(s) ds &<\int_0^{L_{\M}}\frac{2\pi}{L_{\M}}ds,\\
2\pi &<2\pi.
\end{align*}
We get the same result if we assume that $\displaystyle\kappa_{\M}(s)>\frac{2\pi}{L_{\M}}$ for all $s$, therefore $\Sms(\M)$ got at least one cusp singularity.

Let $\M$ be a generic oval.

Let (\ref{ParameterizationM}) be the parameterization of $\M$ and $p(\theta)$ its Minkowski support function.

By Remark \ref{RemKnownWcEq} there are at least $3$ cusps of the Wigner caustic of $\M$, then there exist $\theta_1<\theta_2<\theta_3$ such that $\rho_{\M}(\theta_i)=\rho_{\M}(\theta_i+\pi)$ for $i=1,2,3$ and by Corollary \ref{SingCor} there are at least $4$ cusps of the Constant Width Measure Set of $\M$, then there exists $\varphi_1<\varphi_2$ such that $\displaystyle\rho_{\M}(\varphi_j)+\rho_{\M}(\varphi_j+\pi)=\overline{w}_{\M}$ for $j=1,2$, therefore it is easy to see that there are at least four values of $\theta\in S^1$ such that $\displaystyle\rho_{\M}(\theta)=\frac{1}{2}\overline{w}_{\M}$, hence there are at least $4$ cusps of the Spherical Measure Set of $\M$.
\end{proof}

By Theorem \ref{ThmCuspsSms} we get another proof of the classical Four-Vertex Theorem for ovals.

\begin{cor}\label{Cor4VertexThm}
Let $\M$ be a generic oval. Then $\M$ has at least four vertices.
\end{cor}
\begin{proof}
It is a consequence of Roll's theorem and the fact that there are at least $4$ points on $\M$ in which the curvature is equal to $\displaystyle\frac{2\pi}{L_{\M}}$ -- see Proposition \ref{ThmCurvParallelSms}(i) and Theorem \ref{ThmCuspsSms}.
\end{proof}

\begin{prop}\label{PropDoubleCovSms}
If $\M$ is an oval of constant width, then $\Sms(\M)=\Eq_{\frac{1}{2}}(\M)$. Thus, the map $\M\ni\gamma(\theta)\mapsto\gamma_{\Sms(\M)}(\theta)\in\Sms(\M)$ is the double covering of $\Sms(\M)$.

If $\M$ is a curve with the center of symmetry, then $\Cwms(\M)$ and $\Sms(\M)$ are similar curves and the ratio of symmetry between them is equal to $2$.
\end{prop}
\begin{proof}

Let $p$ be the Minkowski support function of $\M$. Then one can check that the Minkowski support function of $\Sms(\M)$ is in the form:
\begin{align}\label{SupportSms}
p_{\Sms(\M)}(\theta)=p(\theta)-\frac{1}{2}\overline{w}_{\M}.
\end{align}
In terms of coefficients of the Fourier series of $p(\theta)$:
\begin{align}\label{SupportFourierSms}
p_{\Sms(\M)}(\theta)=\sum_{n=1}^{\infty}(a_n\cos n\theta+b_n\sin n\theta).
\end{align}

Then Proposition \ref{PropDoubleCovSms} is a consequence of the support functions of $\Eq_{\frac{1}{2}}(\M)$, $\Cwms(\M)$, $\Sms(\M)$ in terms of coefficients of the Fourier series of $p(
\theta)$ -- see (\ref{SupportFourierWc}), (\ref{SupportFourierCwms}), (\ref{SupportFourierSms}). 
\end{proof}

\begin{cor}\label{CorRotationSms}
Let $\M$ be an oval such that $\Sms(\M)$ is a curve with at most cusp singularities. If $\M$ is not a curve of constant width, then the absolute value of the rotation number of $\Sms(\M)$ is equal to one, otherwise it is equal to $\displaystyle\frac{1}{2}$.
\end{cor}

\begin{prop}
There exists oval $\M$ for which the Constant Width Measure Set has exactly $n$ cusp singularities for $n\geqslant 3$.
\end{prop}
\begin{proof}
Let $n\geqslant 3$ and let 
\begin{align*}
p_{n}(\theta)=\left\{\begin{array}{lll} 
n^2+2+\cos n\theta & \text{ if } & n\equiv 1(\text{mod }2) \\
\left(\frac{n}{2}\right)^2+2+\cos\frac{n\theta}{2} & \text{ if } & n\equiv 0(\text{mod }4) \\
n^2+2+\cos\frac{(n-2)\theta}{2}+\cos\frac{n\theta}{2} & \text{ if } & n\equiv 2(\text{mod }4) \\
\end{array}\right.
\end{align*}
be a Minkowski support function of an oval $\M_{n}$. Then one can check that $\Sms(\M_{n})$ has exactly $n$ cusp singularities.
\end{proof}

\begin{figure}[h]
\centering
\includegraphics[scale=0.4]{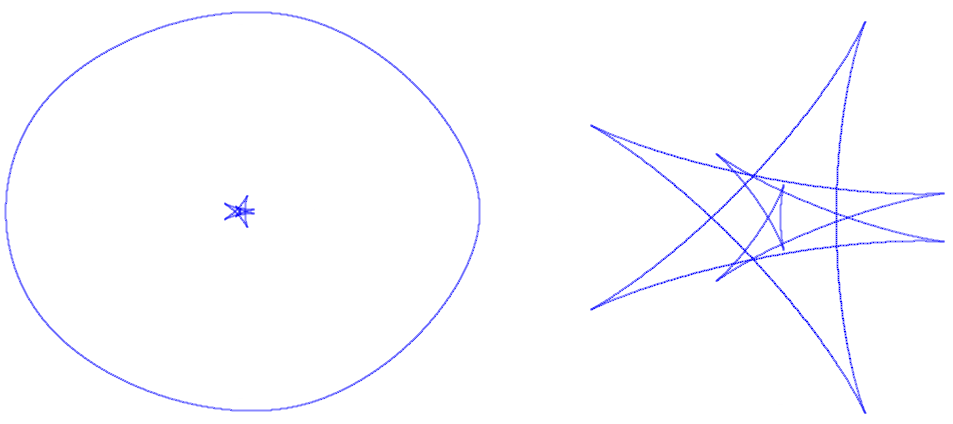}
\caption{An oval $\M_{10}$ and $\Sms(\M_{10})$ which has $10$ cusp singularities. The support function of $M_{10}$ is $p_{10}(\theta)=102+\cos 4\theta+\cos 5\theta$}
\label{PictureM10}
\end{figure}

In Figure \ref{PictureM10} we illustrate example of an oval $\M$ for which $\Sms(\M)$ has exactly $10$ cusp singularities.

\begin{prop}
Let $\M$ be an oval. Then
\begin{align}
L_{\Sms(\M)}\leqslant 2L_{\M}
\end{align}
and if $\M$ is an oval of constant width, then
\begin{align}\label{SmsLenCW}
L_{\Sms(\M)}\leqslant L_{\M}
\end{align}
\end{prop}
\begin{proof}
Let $p(\theta)$ be the Minkowski support function of $\M$ and let $p_{\Sms(\M)}$ (like in (\ref{SupportSms})) be the support function of $\Sms(\M)$.
Let $\M$ is not an oval of constant width. Then
\begin{align}\label{ProofLenSms}
L_{\Sms} &= \int_0^{2\pi}|\gamma'_{\Sms(\M)}(\theta)|\d\theta\\
\nonumber	&=\int_0^{2\pi}|p_{\Sms(\M)}(\theta)+p''_{\Sms(\M)}(\theta)|\d\theta\\
\nonumber	&=\int_0^{2\pi}\left|p(\theta)-\frac{1}{2}\overline{w}_{\M}+p''(\theta)\right|\d\theta\\
\nonumber	&\leqslant\int_0^{2\pi}\left|p(\theta)+p''(\theta)\right|\d\theta+\frac{1}{2}\int_0^{2\pi}\overline{w}_{\M}\d\theta\\
\nonumber	&=L_{\M}+L_{\M}=2L_{\M}.
\end{align}
By (\ref{ProofLenSms}) and Proposition \ref{PropDoubleCovSms} the inequality (\ref{SmsLenCW}) holds.
\end{proof}


\section{Isoperimetric equalities}
\begin{thm}(Isoperimetric equality 1)\label{ThmIsoEq1}
Let $\M$ be a closed regular simple convex planar curve. Then
\begin{align}\label{IsoEq}
L_{\M}^2=4\pi A_M+8\pi\left|\widetilde{A}_{\Eq_{\frac{1}{2}}(\M)}\right|+\pi\left|\widetilde{A}_{\Cwms(\M)}\right|,
\end{align}
where $L_{\M}, A_{\M}, \widetilde{A}_{\Eq_{\frac{1}{2}}(\M)}, \widetilde{A}_{\Cwms(\M)}$ are the length of $\M$, the area bounded by $\M$, the oriented area of the Wigner caustic of $\M$ and the oriented area of the Constant Width Measure Set of $\M$, respectively.
\end{thm}
\begin{proof}
Let $p(\theta)$ be the Minkowski support function of $\M$ and let (\ref{Fourierofp}) be its Fourier series. Then (\ref{SupportFourierCwms}) is the Fourier series of the Minkowski support function $p_{\Cwms(\M)}(\theta)$ of $\Cwms(\M)$.

The oriented area of $\Cwms(\M)$ in terms of coefficients of the Fourier series of $p$ is equal
\begin{align}\label{FourierAreaCwms}
\widetilde{A}_{\Cwms(\M)} &=
\frac{1}{2}\int_0^{2\pi}\Big(p^2_{\Cwms(\M)}(\theta)-p'^2_{\Cwms(\M)}(\theta)\Big)\d\theta\\
\nonumber &= -2\pi\sum_{\substack{n=2,\\ n\text{ is even}}}^{\infty}(n^2-1)(a_n^2+b_n^2).
\end{align}

Then by (\ref{Lengthofmfourier}), (\ref{Areaofmfourier}), (\ref{WCAreaFourier}), (\ref{FourierAreaCwms}) one can easily verify that the equality (\ref{IsoEq}) holds.
\end{proof}

The next theorem is a classical one (see \cite{G4}).

\begin{thm}(Isoperimetric equality 2)\label{ThmIsoEq2}
Let $\M$ be a closed simple regular planar curve. If $\M$ is not a curve of constant width, then
\begin{align}\label{IsoEq2}
L_{\M}^2=4\pi A_M+4\pi\left|\widetilde{A}_{\Sms(\M)}\right|,
\end{align}
otherwise
\begin{align}\label{IsoEq3}
L_{\M}^2=4\pi A_M+8\pi\left|\widetilde{A}_{\Sms(\M)}\right|,
\end{align}
where $L_{\M}, A_{\M}, \widetilde{A}_{\Sms(\M)}$ are the length of $\M$, the area bounded by $\M$ and the oriented area of the Spherical Measure Set of $\M$, respectively.
\end{thm}
\begin{proof}
Let $\M$ be not a curve of constant width.

Let $[0,L_{\M})\ni s\mapsto\gamma(s)\in\M$ be an arc length parameterization of $\M$. Then
$\displaystyle [0,L_{\M})\ni s\mapsto\gamma_{\Sms(\M)}(s)=\gamma(s)+\frac{L_{\M}}{2\pi}\cdot\mathbbm{n}(s)\in\Sms(\M)$ is a parameterization of $\M$.

Let $\mathbbm{t}, \mathbbm{n}, \mathbbm{b}$ be the Frenet frame of $\M$.

By Green's theorem the oriented area of $\Sms(\M)$ is equal to
\begin{align}\label{AreaSms}
\widetilde{A}_{\Sms(\M)} &=\frac{1}{2}\int_0^{L_{\M}}\left(\gamma_{\Sms(\M)}(s)\times\gamma'_{\Sms(\M)}(s)\right)\cdot \mathbbm{b}(s)ds\\
\nonumber&=\frac{1}{2}\int_0^{L_{\M}}\left(\gamma(s)+\frac{L_{\M}}{2\pi}\mathbbm{n}(s)\right)\times\left(1-\frac{L_{\M}}{2\pi}\kappa(s)\right)\mathbbm{t}(s)\cdot \mathbbm{b}(s)ds
\end{align}
\begin{align}
\nonumber&=\frac{1}{2}\int_0^{L_{\M}}\gamma(s)\times \mathbbm{t}(s)\cdot \mathbbm{b}(s)ds-\frac{L_{\M}}{2\pi}\int_0^{L_{\M}}\mathbbm{b}(s)\cdot \mathbbm{b}(s)ds+\\
\nonumber &\ \ \ +\frac{1}{2}\cdot\frac{L_{\M}^2}{4\pi^2}\int_0^{L_{\M}}\kappa(s)ds\\
\nonumber &=A_{\M}-\frac{L_{\M}}{2\pi}\cdot L_{\M}+\pi\cdot\frac{L_{\M}^2}{4\pi^2}=A_{\M}-\frac{L_{\M}^2}{4\pi}.
\end{align}
\end{proof}

By the isoperimetric equalities one can get following corollaries.

\begin{cor}[\cite{Z2}]
A closed regular simple convex planar curve  $\M$ is a curve of constant width if and only if 
\begin{align*}
L_{\M}^2=4\pi A_M+8\pi\left|\widetilde{A}_{\Eq_{\frac{1}{2}}(\M)}\right|,
\end{align*}
where $L_{\M}, A_{\M}, \widetilde{A}_{\Eq_{\frac{1}{2}}(\M)}$ are the length of $\M$, the area bounded by $\M$ and the oriented area of the Wigner caustic of $\M$, respectively.
\end{cor}

\begin{cor}
A closed regular simple convex planar curve  $\M$ is a curve which has the center of symmetry if and only if 
\begin{align*}
L_{\M}^2=4\pi A_M+\pi\left|\widetilde{A}_{\Cwms(\M)}\right|,
\end{align*}
where $L_{\M}, A_{\M}, \widetilde{A}_{\Cwms(\M)}$ are the length of $\M$, the area bounded by $\M$ and the oriented area of the Constant Width Measure Set of $\M$, respectively.
\end{cor}

\begin{cor}\label{CorAreaSmsWcCwms}
Let $\M$ be a closed simple convex planar curve. If $\M$ is not a curve of constant width then
\begin{align}
4\left|\widetilde{A}_{\Sms(\M)}\right| &=8\left|\widetilde{A}_{\Eq_{\frac{1}{2}}(\M)}\right|+\left|\widetilde{A}_{\Cwms(\M)}\right|,
\end{align}
otherwise
\begin{align}
\Sms(\M)=\Eq_{\frac{1}{2}}(\M), \Cwms(\M)=\{\Theta\},
\end{align}
where $\widetilde{A}_{\Eq_{\frac{1}{2}}(\M)}, \widetilde{A}_{\Cwms(\M)}, \widetilde{A}_{\Sms(\M)}$ are the oriented area of the Wigner caustic of $\M$, the oriented area of the Constant Width Measure Set of $\M$ and the oriented area of the Spherical Measure Set of $\M$, respectively.
\end{cor}

By direct calculations we get the following proposition.

\begin{prop}
Let $\M$ be a closed simple convex planar curve. If $\M$ is not a curve of constant width then
\begin{align}
L_{\Sms(\M)}\leqslant\frac{1}{2}L_{\Cwms(\M)}+2L_{\Eq_{\frac{1}{2}}(\M)}.
\end{align}
\end{prop}

\begin{prop}\label{PropOrientationCwmsSms}
Let $\M$ be a generic regular simple closed curve. Then the orientation of $\Sms(\M)$ is reversed against the orientation of $\M$.
If $\M$ is a generic oval, then the orientation of $\Cwms(\M)$ is reversed against the orientation of $\M$.
\end{prop}
\begin{proof}
Let $\M$ be positively oriented. Then the negatively orientation of $\Cwms(\M)$ and $\Sms(\M)$ is a consequence of (\ref{FourierAreaCwms}) and (\ref{AreaSms}).
\end{proof}

\begin{figure}[h]
\centering
\includegraphics[scale=0.15]{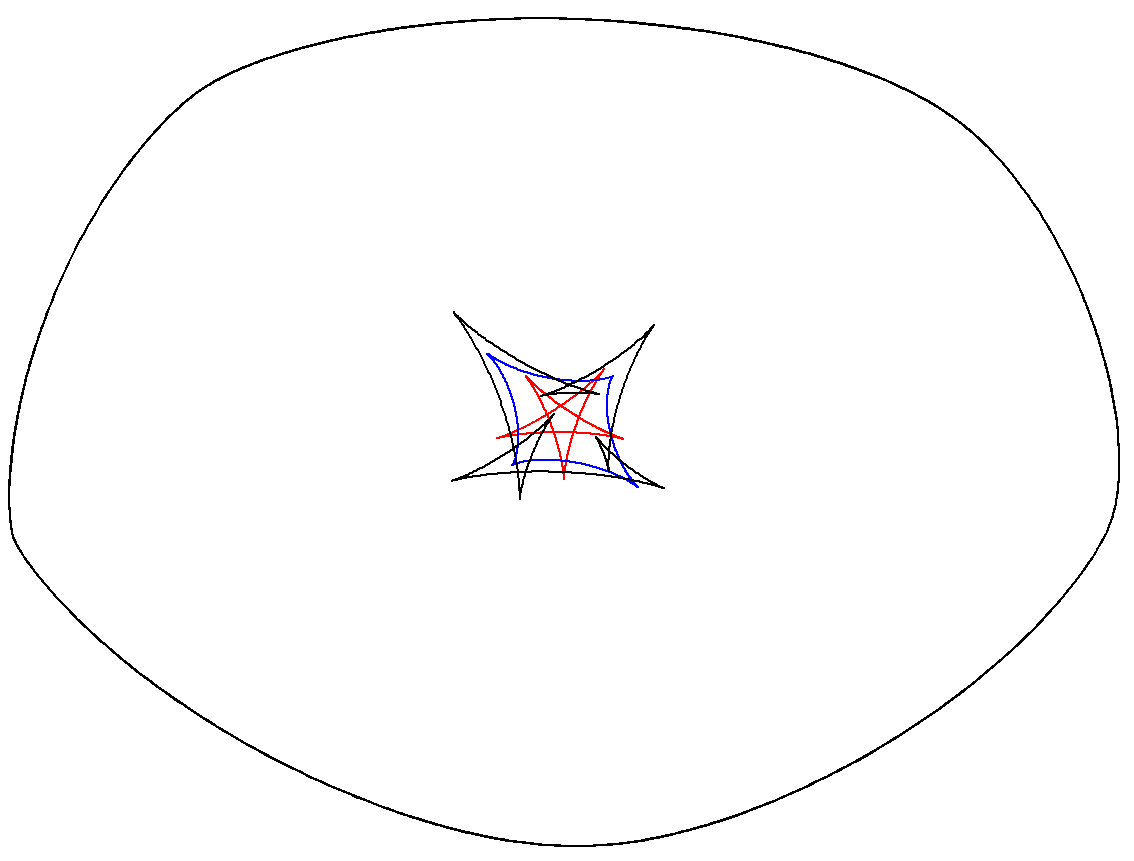}
\includegraphics[scale=0.18]{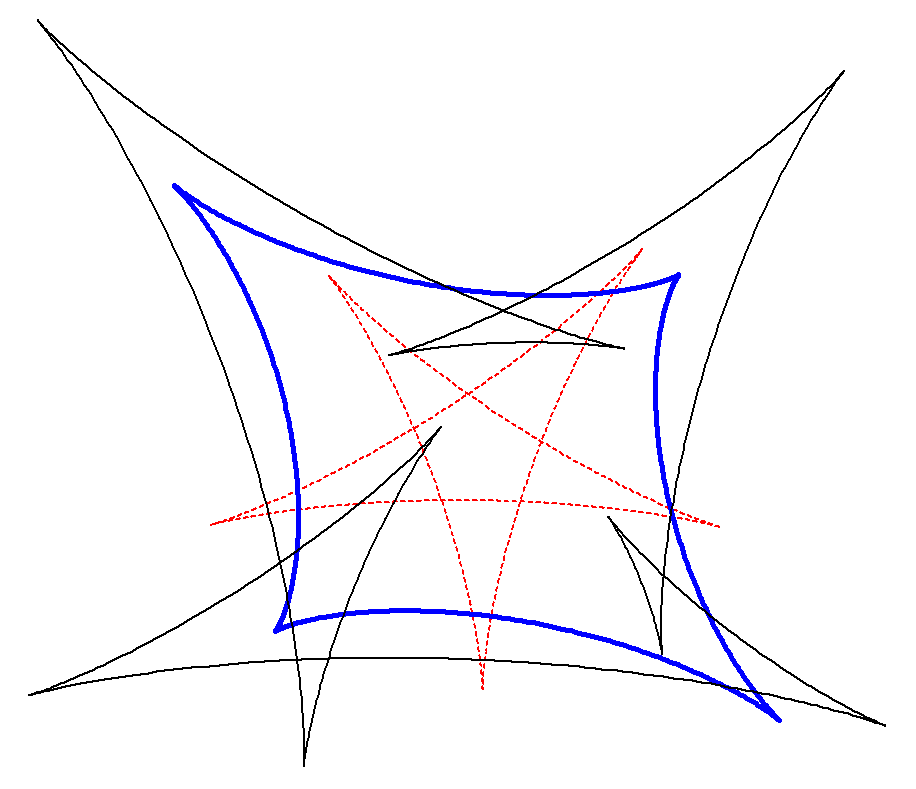}
\caption{A curve $\M$ with the support function \eqref{eq:supAll}, $\Eq_{\frac{1}{2}}(\M)$, $\Cwms(\M)$, and $\Sms(\M)$}
\label{PictureAll}
\end{figure}

Figure \ref{PictureAll} illustrates an example of an oval $\M$ for which 
\begin{align}
\label{eq:supAll}p(\theta)=115+10\cos 2\theta+\frac{1}{3}\cos 3\theta+\sin 4\theta-3\sin 5\theta
\end{align} is its Minkowski support function. In this case $\Eq_{\frac{1}{2}}(\M)$ has $5$ cusp singularities (the dashed line), $\Cwms(\M)$ has $4$ cusp singularities (the bold line) and $\Sms(\M)$ has $10$ cusp singularities. Furthermore the support functions of these curves are as follows:
\begin{align*}
p_{\Eq_{\frac{1}{2}}(\M)}(\theta) &=\frac{1}{3}\cos 3\theta-3\sin 5\theta,\\
p_{\Cwms(\M)}(\theta) &=10\cos 2\theta+\sin 4\theta, \\
p_{\Sms(\M)}(\theta) &=10\cos 2\theta+\frac{1}{3}\cos 3\theta+\sin 4\theta-3\sin 5\theta.
\end{align*}


\section{The stability of improved isoperimetric inequalities}

A \textit{$2$-dimensional convex body} in $\mathbb{R}^2$ is a bounded convex subset in $\mathbb{R}^2$ which is closed and has the non-empty interior. By $\mathcal{C}^2$ we will denote the set of all $n$-dimensional convex bodies.

An inequality in the convex geometry can be written as
\begin{align}\label{IneqConvexGeometry} 
\Phi(K)\geqslant 0,
\end{align}
where $\Phi:\mathcal{C}^2\to\mathbb{R}$ is a function and the inequality (\ref{IneqConvexGeometry}) holds for all $K$ in $\mathcal{C}^2$. Let $\mathcal{C}^2_{\Phi}$ be a subset of $\mathcal{C}^2$ for which the equality in (\ref{IneqConvexGeometry}) holds. 

Let $L_{\partial K}$ (respectively $A_{\partial K}$) denotes the length of the boundary of $K$ (respectively the area enclosed by $\partial K$, i.e. the area of $K$).

In this section we will study stability properties associated with (\ref{IneqConvexGeometry}). We ask if $K$ must be close to a member of $\mathcal{C}^2_{\Phi}$ whenever $\Phi(K)$ is close to zero. If $d:\mathcal{C}^2\times\mathcal{C}^2\to\mathbb{R}$ satisfies two following conditions:
\begin{enumerate}[(i)]
\item $d(K,L)\geqslant 0$ for all $K,L\in\mathcal{C}^2$,
\item $d(K,L)=0$ if and only if $K=L$,
\end{enumerate}
then $d$ denotes in some sense the deviation between two convex bodies.

If $\Phi, \mathcal{C}^2_{\Phi}$ and $d$ are given, then the \textit{stability problem} associated with (\ref{IneqConvexGeometry}) is as follows.

\begin{prob}Find positive constants $c,\alpha$ such that for each $K\in\mathcal{C}^2$, there exists $N\in C^2_{\Phi}$ such that 
\begin{align}\label{StabilityIneq}
\Phi(K)\geqslant cd^{\alpha}(K,N).
\end{align}
\end{prob}

Let $p_{\partial K}$ and $p_{\partial N}$ be support functions of convex bodies $K$ and $N$, respectively. Usually to measure the deviation between $K$ and $N$ one can use the \textit{Hausdorff distance}, $d_{\infty}$, and the measure that corresponds to the $L_2$-metric in the function space, $d_2$, which are given by the following formulas: 
\begin{align}\label{HausdorffDistance}
d_{\infty}(K,N)&=\max_{\theta}\Big|p_{\partial K}(\theta)-p_{\partial N}(\theta)\Big|,\\
\label{LTwoDistance}
d_2(K,N)&=\left(\int_0^{2\pi}\Big|p_{\partial K}(\theta)-p_{\partial N}(\theta)\Big|^2\d\theta\right)^{\frac{1}{2}}.
\end{align}

It is easy to see that $d_{\infty}(K,N)=0$ (or $d_2(K,N)=0$) if and only of $K=N$.

\begin{defn}[\cite{Z2}]
Let $p_{\M}$ be the Minkowski support function of a positively oriented oval $\M$ of length $L_{\M}$. Then 
\begin{align}\label{SupportWM}
p_{W_{\M}}(\theta)=\frac{L_{\M}}{2\pi}+\frac{p_{\M}(\theta)-p_{\M}(\theta+\pi)}{2}
\end{align}
will be the support function of a curve $W_{\M}$ which will be called the \textit{Wigner caustic type curve associated with $\M$}.
\end{defn}

\begin{prop}[\cite{Z2}]
Let $W_{\M}$ be the Wigner caustic type curve associated with an oval $\M$. If $W_{\M}$ is an oval, then
\begin{enumerate}[(i)]
\item $W_{\M}$ is an oval of constant width,
\item $L_{W_{\M}}=L_{\M}$,
\item $\Eq_{\frac{1}{2}}(W_{\M})=\Eq_{\frac{1}{2}}(\M)$,
\item $A_{W_M}\geqslant A_{\M}$ and the equality holds if and only if $\M$ is a curve of constant width,
\item ${W_M}=\M$ if and only if $\M$ is a curve of constant width.
\end{enumerate}
\end{prop}

\begin{rem}\label{RemarkMistake}
Proposition 4.3(i) in \cite{Z2} wrongly states that $W_{\M}$ is always an oval. One can check that the Wigner caustic of an oval $M$ evolving in the perimeter-preserving flow introduced by L. Gao and S. Pan in \cite{GS1} is fixed. Furthermore limiting curve in this flow is exactly $W_{\M}$. By Main Theorem in \cite{GS1} we know that if 
\begin{align}\label{TheCond}
\kappa_{\max}<3\kappa_{\min},
\end{align}
where $\kappa_{\max}:=\max_\theta\kappa_M(\theta)$, $\kappa_{\min}:=\min_\theta\kappa_M(\theta)$, then $W_{\M}$ is an oval.

Example 0 in \cite{GS1} illustrates an oval $M$ with the support function in the form $p(\theta)=10-\cos 2\theta+\cos 3\theta+\frac{1}{8}\cos 5\theta$ for which $W_{\M}$ is a singular hedgehog.
\end{rem}

In \cite{Z2} we show the stability properties of Theorem \ref{ImprovedIsoperimetricIneq1}. By Remark \ref{RemarkMistake} we had to add the condition that $W_{\M}$ is an oval.

\begin{thm}[\cite{Z2}]\label{ThmStabIneqMax1}
Let $K$ be strictly convex domain of area $A_{\partial K}$ and perimeter $L_{\partial K}$ and let $\widetilde{A}_{E_{\frac{1}{2}}(\partial K)}$ denote the oriented area of the Wigner caustic of $\partial K$. Let $W_{K}$ denote the convex body for which $\partial W_{K}$ is the Wigner caustic type curve associated with $\partial K$. If $W_M$ is an oval, then
\begin{align}\label{StabIneqMax}
L_{\partial K}^2-4\pi A_{\partial K}-8\left|\widetilde{A}_{\Eq_{\frac{1}{2}}(\partial K)}\right|\geqslant 4\pi^2 d_{\infty}^2(K, W_K),
\end{align}
where equality holds if and only if $\partial K$ is a curve of constant width.
\end{thm}

\begin{thm}[\cite{Z2}]\label{ThmStabIneqMax2}
Under the same assumptions of Theorem \ref{ThmStabIneqMax1}, one gets
\begin{align}\label{StabIneqMax}
L_{\partial K}^2-4\pi A_{\partial K}-8\pi\left|\widetilde{A}_{\Eq_{\frac{1}{2}}(\partial K)}\right|\geqslant 6\pi d_2^2(K, W_K),
\end{align}
where equality holds if and only if $\partial K$ is a curve of constant width, or the Minkowski support function of $\partial K$ is in the form 
\begin{align*} p_{\partial K}(\theta)=a_0+a_2\cos 2\theta+b_2\sin 2\theta+\sum_{\substack{n=1,\\ n\text{ is odd}}}^{\infty}(a_{n}\cos n\theta+b_{n}\sin n\theta).
\end{align*}
\end{thm}

By Theorem \ref{ThmStabIneqMax1}, Theorem \ref{ThmStabIneqMax2} and Theorem \ref{ThmIsoEq1} we can get the following result.

\begin{cor}\label{CorBoundCwms}
Under the same assumptions of Theorem \ref{ThmStabIneqMax1}, one gets
\begin{align}
\left|\widetilde{A}_{\Cwms(\partial K)}\right| &\geqslant \max\Big\{4\pi d_{\infty}^2(K, W_K),\ 6 d_2^2(K, W_K)\Big\}.
\end{align}
\end{cor}

By Theorem \ref{ThmIsoEq1} we can easily obtain the second improved isoperimetric inequality (see Theorem \ref{ThmImpIsoIneq2}), study its stability property and obtain similar results like these in Corollary \ref{CorBoundCwms}, but for the Wigner caustic of an oval.

\begin{thm}(Improved isoperimetric inequality 2)\label{ThmImpIsoIneq2}
Let $\M$ be a closed regular simple convex planar curve. Then
\begin{align}\label{IsoEq}
L_{\M}^2\geqslant 4\pi A_M+\pi\left|\widetilde{A}_{\Cwms(\M)}\right|,
\end{align}
where $L_{\M}, A_{\M}, \widetilde{A}_{\Cwms(\M)}$ are the length of $\M$, the area bounded by $\M$ and the oriented area of the Constant Width Measure Set of $\M$, respectively.
The equality in (\ref{IsoEq}) holds if and only if $\M$ is a centrally symmetric.
\end{thm}

From now one let us assume that $n=2$ and by Theorem \ref{ThmImpIsoIneq2} let
\begin{align}\label{IneqConvexProblem2}
\Phi(K)=L_{\partial K}^2-4\pi A_{\partial K}-\pi\left|\widetilde{A}_{\Cwms(\partial K)}\right|.
\end{align}
By Theorem \ref{ThmImpIsoIneq2} one can see that $C^2_{\Phi}$ consists of centrally symmetric bodies.

Let $K$ be the convex body in $\mathbb{R}^2$. Let $\vec u$ be a unit vector in $\mathbb{R}^2$. Then the \textit{Steiner point} which can be defined in terms of the Minkowski support function of $\partial K$ is as follows:
\begin{align}
\vec s(K)=\frac{1}{\pi}\int_0^{2\pi}\vec u(\theta)\cdot p_{\partial K}(\theta)\d\theta.
\end{align}

In terms of coefficients of the Fourier series of $p_{\partial K}$ one can notice that $\vec s(K)=(a_1, b_1)$.

For more details on the Steiner point of a convex body see \cite{G5, S2, S3}.

\begin{defn}
Let $K$ be the convex body. Then $\displaystyle S_{K}=\frac{1}{2}(K+(-K))$ is centrally symmetric convex body, where $+$ is the Minkowski addition and \linebreak $-K=\{-k\ |\ k\in K\}$. Often $S_{K}$ is called the \textit{Steiner symmetral of K}, and the process of generating the set from $K$ is known as \textit{symmetrization}. Let as translate $S_{K}$ such that the center of $S_{K}$ becomes the Steiner point of $K$.
\end{defn}

The support function of $\partial S_{K}$ is equal to
\begin{align}\label{SupportSM}
p_{\partial S_{K}}(\theta)=\vec s(K)+\frac{p_{\partial K}(\theta)+p_{\partial K}(\theta+\pi)}{2}
\end{align}
and in terms of coeeficients of the Fourier series of $p_{\partial K}$:
\begin{align}\label{SupportSMFourier}
p_{\partial S_{K}}(\theta)=a_0+a_1\cos\theta+b_1\sin\theta+\sum_{\substack{n=2,\\ n\text{ is even}}}^{\infty}\left(a_n\cos n\theta+b_n\sin n\theta\right).
\end{align}

\begin{thm}\label{ThmStabIneqMax2}
Let $K$ be strictly convex domain of the area $A_{\partial K}$ and the perimeter $L_{\partial K}$ and let $\widetilde{A}_{\Cwms(\partial K)}$ denote the oriented area of the Constant Width Measure Set of $\partial K$. Let $S_{K}$ be the Steiner symmetral of $K$. Then
\begin{align}\label{StabIneqMaxCwms}
L_{\partial K}^2-4\pi A_{\partial K}-\pi\left|\widetilde{A}_{\Cwms(\partial K)}\right|\geqslant 8\pi^2 d_{\infty}^2(K, S_K),
\end{align}
where equality holds if and only if $\partial K$ is centrally symmetric.
\end{thm}

\begin{proof}

By (\ref{Fourierofp}), (\ref{Fourierofpprime}), (\ref{CauchyFormula}), (\ref{BlaschkeFormula}),  one can get the Fourier series of $\Phi$ (see (\ref{IneqConvexProblem2})):

\begin{align}\label{FourierOfPhi}
\Phi(K) &=L_{\partial K}^2-4\pi A_{\partial K}-8\pi\left|\widetilde{A}_{E_{\frac{1}{2}}(\M)}\right|=2\pi^2\sum_{\substack{n=3,\\ n\textit{ is odd}}}^{\infty}(n^2-1)(a_n^2+b_n^2).
\end{align}

One can check that $|a_n\cos n\theta+b_n\sin n\theta|\leqslant\sqrt{a_n^2+b_n^2}$ and then by (\ref{HausdorffDistance}) and H\"older's inequality:

\begin{align*}
d_{\infty}(K, S_K) &= \max_{\theta}\Big|p_{\partial K}(\theta)-p_{\partial S_K}(\theta)\Big|   =\max_{\theta}\left|\sum_{\substack{n=3,\\ n\textit{ is odd}}}^{\infty}(a_n\cos n\theta+b_n\sin n\theta)\right|\\
	&\leqslant\max_{\theta}\left(\sum_{\substack{n=3,\\ n\textit{ is odd}}}^{\infty}|a_n\cos n\theta+b_n\sin n\theta|\right)\\
	&\leqslant\sum_{\substack{n=3,\\ n\textit{ is odd}}}^{\infty}\frac{1}{\sqrt{n^2-1}}\cdot \sqrt{n^2-1}\sqrt{a_n^2+b_n^2}\\
	&\leqslant\sqrt{\sum_{\substack{n=3,\\ n\textit{ is odd}}}^{\infty}\frac{1}{n^2-1}}\cdot \sqrt{\sum_{\substack{n=3,\\ n\textit{ is odd}}}^{\infty}(n^2-1)(a_n^2+b_n^2)}=\sqrt{\frac{1}{4}}\cdot\sqrt{\frac{\Phi(K)}{2\pi^2}}.
\end{align*}
And the equality holds if and only if $a_{2m+1}=b_{2m+1}=0$ for all $m\in\mathbb{N}$, so $\partial K$ is centrally symmetric.
\end{proof}

\begin{thm}\label{ThmStabIneqMax3}
Under the same assumptions of Theorem \ref{ThmStabIneqMax2}, one gets
\begin{align}\label{StabIneqMax}
L_{\partial K}^2-4\pi A_{\partial K}-8\left|\widetilde{A}_{\Cwms(\partial K)}\right|\geqslant 16\pi d_2^2(K, S_K),
\end{align}
where equality holds if and only if $\partial K$ is centrally symmetric or the Minkowski support function of $\partial K$ is in the form 
\begin{align*} 
p_{\partial K}(\theta) =&\ a_0+a_1\cos\theta+b_1\sin\theta+a_3\cos 3\theta+b_3\sin 3\theta\\
	&+\sum_{\substack{n=2,\\ n\text{ is even}}}^{\infty}(a_{n}\cos n\theta+b_{n}\sin n\theta).
\end{align*}
\end{thm}

\begin{proof}
By (\ref{SupportSMFourier}) and (\ref{FourierOfPhi})
\begin{align*}
d_2^2(K,S_K) &=\int_0^{2\pi}\Big|p_{\partial K}(\theta)-p_{\partial S_K}(\theta)\Big|^2\d\theta=\int_0^{2\pi}\left|\sum_{\substack{n=3,\\ n\textit{ is odd}}}^{\infty}(a_n\cos n\theta+b_n\sin\theta)\right|^2\d\theta\\
	&=\pi\sum_{\substack{n=3,\\ n\textit{ is odd}}}^{\infty}(a_n^2+b_n^2)\leqslant\frac{1}{16\pi}\cdot 2\pi^2\sum_{\substack{n=3,\\ n\textit{ is odd}}}^{\infty}(n^2-1)(a_n^2+b_n^2)=\frac{1}{16\pi}\Phi(K).
\end{align*}
And the equality holds if and only if $a_{2m+1}=b_{2m+1}=0$ for all $m\in\mathbb{N}$, so $\partial K$ is a centrally symmetric curve, or $\displaystyle p_{\partial K}(\theta)=a_0+a_1\cos\theta+b_1\sin\theta+a_3\cos 3\theta+b_3\sin 3\theta+\sum_{\substack{n=2,\\ n\text{ is even}}}^{\infty}(a_{n}\cos n\theta+b_{n}\sin n\theta)$.
\end{proof}

\begin{cor}\label{CorBoundWc}
Under the same assumptions of Theorem \ref{ThmStabIneqMax2}, one gets
\begin{align}
\label{CorBoundWc1}\left|\widetilde{A}_{\Eq_{\frac{1}{2}}(\partial K)}\right| \geqslant \max &\left\{\pi d_{\infty}^2(K, S_K),\ 2d_2^2(K, S_K)\right\},
\end{align}
\begin{align}
\label{CorBoundSms1} \left|\widetilde{A}_{\Sms(\partial K)}\right| \geqslant \max &\left\{
2\pi d_{\infty}^2(K, S_K)+\pi d_{\infty}^2(K, W_K),\ 2\pi d_{\infty}^2(K, S_K)+\frac{3}{2} d_2^2(K, W_K),\right.\\ \nonumber& \left. 4 d_2^2(K, S_K)+\pi d_{\infty}^2(K, W_K),\ 4 d_2^2(K, S_K)+\frac{3}{2} d_2^2(K, W_K)\right\},
\end{align}
where (\ref{CorBoundSms1}) holds whenever $\partial K$ is not a curve of constant width and if $W_{K}$ is an oval.
\end{cor}
\begin{proof}
It is an easy consequence of isoperimetric equalities (see Theorem \ref{ThmIsoEq1}--\ref{ThmIsoEq2} and Corollary \ref{CorAreaSmsWcCwms}) and stability properties of improved isoperimetric inequalities (see Corollary \ref{CorBoundCwms} and Theorem \ref{ThmStabIneqMax2}--\ref{ThmStabIneqMax3}).
\end{proof}

\section*{Acknowledgements}
The author is sincerely grateful to Professor Wojciech Domitrz for the valuable discussions.

\bibliographystyle{amsalpha}

\end{document}